\documentclass[reqno]{amsart} 
\usepackage{amsfonts, amsmath, amsthm, amssymb, latexsym, graphicx}

\theoremstyle{plain}
\newtheorem{theorem}{Theorem}[section]
\newtheorem{corollary}[theorem]{Corollary}
\newtheorem{lemma}[theorem]{Lemma}
\newtheorem{proposition}[theorem]{Proposition}
\theoremstyle{definition}

\theoremstyle{remark}

\numberwithin{equation}{section}

\title[]{An estimate for $\beta$-Hermite ensembles via the zeros of Hermite polynomials}

\author{Michael Voit}
\address{Fakult\"at Mathematik, Technische Universit\"at Dortmund,
          Vogelpothsweg 87,
          D-44221 Dortmund, Germany}
\email{michael.voit@math.tu-dortmund.de}


\begin{document}
\date{\today}

\begin{abstract}
  Let $X$ be an $N$-dimensional random vector which describes the ordered eigenvalues
  of a  $\beta$-Hermite  ensemble, and let  $z$
  the vector containing the ordered zeros of the Hermite poynomial $H_N$.
  We  present an explicit estimate for $P(\|X-z\|_2\ge\epsilon)$ for small $\epsilon>0$ and large parameters $\beta$.
  The proof is based on  a central limit theorem for these ensembles for $\beta\to\infty$
  with explicit eigenvalues of the covariance matrices of  the limit. The estimate is similar to previous estimates of Dette and Imhof (2009).
\end{abstract}

\maketitle

\section{Introduction} 

For large parameters, the three classical $\beta$-random matrix ensembles  of dimension $N\ge2$ with ordered entries
are related with 
deterministic vectors which contain the ordered zeros of some associated classical orthogonal poynomial of order $N$
as follows:

\smallskip
\rm{(1)} {\bf Hermite ensembles:} Let  $\beta=2k>0$ be a constant (both parameters are used in the the literature;
see e.g. \cite{AKM1, AV, D, DI, DE1, DE2, F, GK, Me, V}).
Then the associated  Hermite ensemble may be regarded as a  random vector $X_{k,N}$
with values in the closed Weyl chamber
$C_N^A:=\{x\in \mathbb R^N: \quad x_1\ge x_2\ge\ldots\ge x_N\}$ with density
\begin{equation}\label{density-A}
c_k^A   e^{-\|x\|_2^2/2} \prod_{1\le i<j\le N}(x_i-x_j)^{2k}
\end{equation}
with the well-known normalization
\begin{equation}\label{const-A}
 c_k^A:= \Bigl(\int_{C_N^A}  e^{-\|y\|_2^2/2}\cdot \prod_{i<j} (y_i-y_j)^{2k} \> dy\Bigr)^{-1}
 =\frac{N!}{(2\pi)^{N/2}} \cdot\prod_{j=1}^{N}\frac{\Gamma(1+k)}{\Gamma(1+jk)}.
\end{equation}
The aim of this note is to compare $X_{k,N}$ with the deterministic vector $\sqrt{2k}\cdot {\bf z}_N\in C_N^A$, where
the entries of $ {\bf z}_N$  are
 the ordered zeros of the  Hermite polynomial $H_N$ 
where the  $(H_N)_{N\ge 0}$ are orthogonal w.r.t.
the density  $e^{-x^2}$ as e.g.~ in \cite{S}.

\smallskip
\rm{(2)}  {\bf Laguerre ensembles:} Let $k_1,k_2>0$ be constants and $k=(k_1,k_2)$. Define
 the associated  Laguerre ensemble as random vector $X_{k,N}$ with values in 
 $$C_N^B:=\{x\in \mathbb R^N: \quad x_1\ge x_2\ge\ldots\ge x_N\ge0\}$$
with density
\begin{equation}\label{density-B}
c_k^B   e^{-\|x\|_2^2/2}\cdot \prod_{i<j}(x_i^2-x_j^2)^{2k_2}\cdot \prod_{i=1}^N x_i^{2k_1}\end{equation}
with some well-known normalization $c_k^B>0$.
 We now write
 $k=(k_1,k_2)=(\nu\cdot\beta,\beta)$ with $\nu>0$ fixed and   and  $\beta\to\infty$.
Then the results in \cite{AKM2, AHV, AV, DE2, DI,  V}  motivate to compare   $X_{k,N}$ with the  vector
$\sqrt{\beta}\cdot {\bf r}_N\in C_N^B$ where
$${\bf r}_N=(\sqrt z_1, \ldots,\sqrt z_N)$$
for 
the ordered zeros $ z_i$ of the Laguerre polynomial   $L_N^{(\nu-1)}$, and  where
the  $L_N^{(\nu-1)}$ are orthogonal  w.r.t. the density $e^{-x}\cdot x^{\nu-1}$ on $]0,\infty[$  for $\nu>0$.

\smallskip
\rm{(3)} {\bf Jacobi ensembles:}
For the parameters $\kappa>0$ and $\alpha,\beta >-1$, the $\beta$-Jacobi ensembles can be regarded as
random vectors $X_{\alpha,\beta,\kappa;N}$ with values in the alcove $A_N:=\{x\in\mathbb R^N: \> -1\le x_1\le\cdots\le x_N\le 1\}$ with densities
\begin{equation}\label{densitiybetaJacobi}
 	C_{\kappa,\alpha,\beta, N,J}
\prod_{1\leq i<j\leq N}(x_j-x_i)^{\kappa}\prod_{i=1}^N(1-x_i)^{\frac{(\alpha+1) \kappa}{2}-\frac12}(1+x_i)^{\frac{(\beta+1)\kappa}{2}-\frac12}\> dx
\end{equation}
 with some known normalizations $C_{\kappa,\alpha,\beta, N,J}B>0$; see \cite{ F,  K, KN, L}.
For fixed $\alpha,\beta$ and $\kappa\to\infty$, the results in \cite{HV, AHV}  motivate to compare $X_{\alpha,\beta,\kappa;N}$  with the  vector
 consisting of the ordered zeros of the Jacobi polynomial $P_N^{(\alpha,\beta)}$ on $[-1,1]$.
\smallskip

For these 3 ensembles there are freezing central limit theorems (FCLTs)
where the covariance matrices of the limits as well as their inverses  and the
associated eigenvalues and eigenvectors are known in terms
of the associated classical orthogonal polynomials; see \cite{DE2, V, AV, AHV, GK, HV}.
These FCLTs can be derived either from the associated tridiagonal models of Dumitriu and Edelman \cite{DE1} as in \cite{DE2}
in the first two cases or in a direct way via Taylor expansions of the densities  \eqref{density-A},
\eqref{density-B}, and \eqref{densitiybetaJacobi}. The approach via the tridiagonal models was used in \cite{DI}
to derive estimates for  $P(\|X_{k,N}-z_N\|_2\ge\epsilon)$ for the ensembles $X_{k,N}$ and corresponding vectors ${\bf z}_N$
above (where $k$ stands for the respective parameters).
In this note we use the  approach  via the eigenvalues of the covariance matrices of the limits and a
refined Taylor expansion
in order to derive estimates for  $P(\|X_{k,N}-z_N\|_2\ge\epsilon)$ in the Hermite case; see Proposition \ref{main-result-a}
and Corollary \ref{main-corr-a}  for details.
It turns out that, roughly, our results are better than those in \cite{DI}
when $k$ is large in comparison to $N$, and worse in the converse case.

The details of the proof Proposition \ref{main-result-a} are quite technical. We point out that however, the approach
can be in principle easily adapted to the Laguere and Jacobi case, where the details of the elementary estimates become even more involved.

\section{Estimates in the Hermite case}

In this Section we derive estimates
for $P\bigl(\|\frac{X_{k,N}}{\sqrt{2k}}-{\bf z}_N\|_2>\epsilon\bigr)$ for Hermite ensembles.
We start with some notations and known results.
We recapitulate from the introduction that for  positive integers $N$, the vector
${\bf z}_N=(z_{1,N},\ldots, z_{N,N})\in C_N^A$  consists of the $N$ ordered zeros of the Hermite polynomial
$H_N$. We need:

\begin{lemma}\label{char-zero-A}
  For a vector $z=(z_1,\ldots,z_N)\in C_N^A$, we have $z= {\bf z}_N$ if and only if
  $$z_i= \sum_{j: j\ne i} \frac{1}{z_i-z_j} \quad\text{for}\quad i=1,\ldots,N.$$
  Moreover, ${\bf z}_N$ satisfies
\begin{equation}\label{potential-z-A}
 -\|{\bf z}_N\|_2^2+ 2\sum_{i<j} \ln(z_{i,N}-z_{j,N})= -\frac{N(N-1)}{2}(1+\ln 2)+ \sum_{j=1}^N j\ln j.
\end{equation}
\end{lemma}

\begin{proof} For the first  statement 
 see Section 6.7 of \cite{S}, or \cite{AKM1} or \cite{AV}.
For (\ref{potential-z-A}) we
 refer to Appendix D and the comments between Eqs.~(58) and (59) in \cite{AKM1}.
\end{proof}

Lemma \ref{char-zero-A} was a main ingredient in \cite{V} for the following central limit theorem for
the Hermite ensembles  for $k\to\infty$. It was also derived  in a  different way
(and with the limit covariance matrix in a different form)
by Dumitriu and Edelman \cite{DE1, DE2} via their tridiagonal $\beta$-ensemble models.
The statement about the  eigenvalues of the covariance matrix in the following CLT can be found in \cite{AV}:

\begin{theorem}\label{clt-main-a}
For each $N\ge 2$,
$X_{k,N} -  \sqrt{2k}\cdot {\bf z}_N$
converges for $k\to\infty$ to the  $N$-dimensional centered  normal distribution $N(0,\Sigma_N)$
with the regular covariance matrix $\Sigma_N$ with $\Sigma_N^{-1}=S_N=(s_{i,j})_{i,j=1}^N$ and
\begin{equation}\label{covariance-matrix-A}
s_{i,j}:=\left\{ \begin{array}{r@{\quad\quad}l}  1+\sum_{l\ne i} (z_{i,N}-z_{l,N})^{-2} & \text{for}\quad i=j \\
   -(z_{i,N}-z_{j,N})^{-2} & \text{for}\quad i\ne j  \end{array}  \right.  . 
\end{equation}
The matrix $S_N$ has the eigenvalues $1,2,\ldots, N$, i.e., the covariance matrix $\Sigma_N$ has the eigenvalues $1,1/2,\ldots, 1/N$.
\end{theorem}

The statement about the eigenvalues in Theorem \ref{clt-main-a} leads to the following:

\begin{corollary}\label{corr-clt-main-a}
  \begin{equation}\label{sum-inverse-herm-1}
      \sum_{i,j: \> i\ne j} (z_{i,N}-z_{j,N})^{-2} =N(N-1)/2\end{equation}
  and
  \begin{equation}\label{sum-inverse-herm-2}
   \sum_{i,j: \> i\ne j} (z_{i,N}-z_{j,N})^{-4} \le \frac{N(N-1)(2N-1)}{12}\le N^3/6.\end{equation}
\end{corollary}

\begin{proof}
(\ref{sum-inverse-herm-1}) follows from
 $$N+\sum_{i,j: \> i\ne j} (z_{i,N}-z_{j,N})^{-2}= tr(S_N)=1+2+\ldots+N=N(N+1)/2 $$
  and can be also found in the literature like in \cite {AKM1}.
For the proof of (\ref{sum-inverse-herm-2})  observe that by (\ref{covariance-matrix-A}) and (\ref{sum-inverse-herm-1}),
  \begin{align} \frac{N(N+1)(2N+1)}{6}&=1^2+2^2+\ldots+N^2 =tr(S_N^2) \notag\\
    &=\sum_{i=1}^N\Bigl(1+ \sum_{j: \> i\ne j}(z_{i,N}-z_{j,N})^{-2} \Bigr)^2 +
    \sum_{i,j: \> i\ne j} (z_{i,N}-z_{j,N})^{-4} \notag\\
 &=2\sum_{i,j: \> i\ne j} (z_{i,N}-z_{j,N})^{-2} + N +2\sum_{i,j: \> i\ne j} (z_{i,N}-z_{j,N})^{-4}  \notag\\
&\quad\quad+ \sum_{i=1}^N\sum_{l,j: \> l\ne j, l\ne i, j\ne i}  (z_{i,N}-z_{j,N})^{-2}(z_{i,N}-z_{l,N})^{-2} \notag\\
 &\ge N(N-1)+N +2\sum_{i,j: \> i\ne j} (z_{i,N}-z_{j,N})^{-4} . \notag
  \end{align}
  This and a short computation now lead to  (\ref{sum-inverse-herm-2}).
\end{proof}

The statement about the eigenvalues in Theorem \ref{clt-main-a} and the method of the proof of Corollary 
\ref{corr-clt-main-a} can be also used to derive lower bounds for
$$M_N:=\min_{1\le i \le  N-1}(z_{i,N}-z_{i+1,N})$$
of order 
$O(1/\sqrt N)$; see \cite{V2}. However, the constant in this approach is slighty worse than known constants
 obtained by other methods. To our knowledge, the best estimate for $M_N$ is the following one in \cite{Kr}:

\begin{lemma}\label{min-distance-a} For all $N\ge2$, $M_N\ge  2/\sqrt N$.
  \end{lemma}

With these results we now derive  the following estimate:

\begin{proposition}\label{main-result-a}
  Let $N\ge2$, $k\ge1$, and $\epsilon>0$ with
  \begin{equation}\label{main-condition-prp-a} \sqrt{\frac{  1+\ln N}{2k}}  \le \epsilon\le \frac{1}{2\sqrt N}  .\end{equation}
  Then
\begin{equation}\label{main-estimate-prp-a}
 P\Bigl(\|\frac{X_{k,N}}{\sqrt{2k}}-{\bf z}_N\|_2>\epsilon\Bigr)
 \le \frac{32}{3}\epsilon^4kN^3-\frac{N-1}{26\cdot k}
 + E\cdot \frac{\sqrt e\sqrt N}{2}
(2k\epsilon^2+1)e^{-k\epsilon^2}
 \end{equation}
with
$$E:=exp\Bigl(\frac{1}{k}\Bigl(\frac{N-1}{26}- \frac{32}{3} \epsilon^4k^2N^3\Bigr)\Bigr).$$
\end{proposition}

Clearly, this result is interesting only for $k\ge 2N(1+\ln N)$.

The terms in the r.h.s. of (\ref{main-estimate-prp-a}) have the following interpretation: By the proof below,
  the term $\frac{32}{3}\epsilon^4kN^3-\frac{N-1}{26\cdot k}$ is technical and caused by the quality of the comparison
  of the densities in the CLT \ref{clt-main-a}; the factor $N^3$ there is slightly annoying,
  but it cannot be impoved in an essential way by the very approach. 
  Moreover, $E$
  is close to 1 in all interesting cases and thus not relevant, while $\frac{\sqrt e\sqrt N}{2}
  (2k\epsilon^2+1)e^{-k\epsilon^2}$ is caused by the form of the (inverse) covariance matrix  in the CLT \ref{clt-main-a}
  and cannot be improved in an essential way.  

  In order to see Proposition \ref{main-result-a} at work,  consider 
  $\epsilon:= c\Bigl(\frac{\ln k}{k}\Bigr)^{1/2}$ with some $c>0$. Here, condition (\ref{main-condition-prp-a}) means that
  \begin{equation}\label{main-condition-prp-a-mod} \sqrt{\frac{  1+\ln N}{2\ln k}}  \le c\le \frac{1}{2}\sqrt{\frac{k}{N\ln k}}
    .\end{equation}
For $N\ge 2$ and $k\ge e$ this implies that $E\le1$, and (\ref{main-estimate-prp-a}) leads to

\begin{corollary}\label{main-corr-a} Let  $N\ge 2$ and $k\ge e$ and $c>0$ with (\ref{main-condition-prp-a-mod}). Then
  \begin{equation}\label{main-estimate-prp-a-mod}
 P\Bigl(\|\frac{X_{k,N}}{\sqrt{2k}}-{\bf z}_N\|_2>c\Bigl(\frac{\ln k}{k}\Bigr)^{1/2}  \Bigr)
 \le \frac{32}{3}c^4N^3\frac{(\ln k)^{2}}{k}+ 
 \frac{\sqrt e N}{2}
(2c^2 \ln k +1)k^{-c^2}.
 \end{equation}
In particular, for $c=1$,  $k\to\infty$, and uniformly in $N$,
\begin{equation}\label{main-estimate-prp-a-mod2}
P\Bigl(\|\frac{X_{k,N}}{\sqrt{2k}}-{\bf z}_N\|_2>\Bigl(\frac{\ln k}{k}\Bigr)^{1/2}  \Bigr)
= O\Bigl(N^3\frac{(\ln k)^2}{k}\Bigr),  \end{equation}
and for $0<c<1$,
\begin{equation}\label{main-estimate-prp-a-mod3}
P\Bigl(\|\frac{X_{k,N}}{\sqrt{2k}}-{\bf z}_N\|_2>c\Bigl(\frac{\ln k}{k}\Bigr)^{1/2}  \Bigr)
= O\Bigl(N^3\frac{\ln k}{k^{c^2}}\Bigr),  \end{equation}
\end{corollary}

We briefly compare the results above with the following closely related result of  Dette and Imhof \cite{DI} which was obtained
via the tridiagonal  matrix model of  Dumitriu and Edelman  \cite{DE1} for $\beta$-Hermite ensembles:
\begin{equation}\label{dette}
  P\Bigl(\|X_{k,N}-\sqrt{2k}\cdot{\bf z}_N\|_\infty> \epsilon)\le 4Ne^{-\epsilon^2/18}. \end{equation}
This yields for $0<c\le1$
\begin{equation}\label{dette2}
  P\Bigl(\|\frac{X_{k,N}}{\sqrt{2k}}-{\bf z}_N\|_2>\Bigl(\frac{\ln k}{k}\Bigr)^{1/2}  \Bigr)
  = O(N\cdot k^{-c^2/9})\end{equation}
Hence, (\ref{main-estimate-prp-a-mod2}) is worse than  (\ref{dette2}) when $N$ is large compared to $k$,
and better
in the converse case.

In view of these results,
we conjecture that in (\ref{main-estimate-prp-a-mod2}) and  (\ref{main-estimate-prp-a-mod3}), the orders 
$O\Bigl(N\frac{(\ln k)^2}{k}\Bigr)$ for $c=1$ and $O\Bigl(N\frac{\ln k}{k^{c^2}}\Bigr)$ for $0<c<1$ respectively are available.

 We now turn to the proof of the Proposition.

\begin{proof}[Proof of Proposition \ref{main-result-a}]
The proof  uses the same arguments as  the proof of the CLT \ref{clt-main-a} in \cite{V}
at the beginning.
In fact, (\ref{density-A}) yields that
$X_{k,N} -  \sqrt{2k}\cdot {\bf z}_N$ has the Lebesgue density
\begin{align}\label{a-density-detail}
f_k^A(y):=&c_k^A  \cdot exp\Bigl(-\|y +\sqrt{2k}\cdot{\bf z}_N\|_2^2/2+2k \sum_{i<j}\ln\bigl(y_i-y_j+\sqrt{2k}(z_{i,N}-z_{j,N})\bigl) \Bigr)
 \notag\\
=c_k^A &\cdot exp\Bigl(-\|y\|_2^2/2-\sqrt{2k} \langle y,{\bf z}_N\rangle- k\|{\bf z}_N\|_2^2
+2k \sum_{i<j}\ln(\sqrt{2k}(z_{i,N}-z_{j,N})) \Bigr)\times\notag\\
&\times  exp\Bigl(
2k \sum_{i<j}\ln\bigl(1+ \frac{y_i-y_j}{\sqrt{2k}(z_{i,N}-z_{j,N})}\bigr)\Bigr) 
\end{align}
on the shifted cone $C_N^A-\sqrt{2k}\cdot{\bf z}_N$ with $f_k^A(y)=0$ otherwise on $\mathbb R^N$.
In the first case, we write $f_k^A$ as
$$f_k^A(y)= \tilde c_k\cdot h_k(y)$$
with
\begin{equation}\label{def-density-a}
    h_k(y):=  exp\Bigl(-\|y\|_2^2/2 -\sqrt{2k} \langle y,{\bf z}_N\rangle +
    2k \sum_{i<j}\ln\bigl(1+ \frac{y_i-y_j}{\sqrt{2k}(z_{i,N}-z_{j,N})}\bigr)\Bigr)
 \end{equation}   
and  the constant
\begin{equation}\label{a-density-constants-a}
  \tilde c_k:= c_k^A exp\Bigl(-k\|{\bf z}_N \|_2^2 + 2k \sum_{i<j}\ln(\sqrt{2k}(z_{i,N}-z_{j,N}))\Bigr).
\end{equation}   
Using  (\ref{potential-z-A})  and  (\ref{const-A}), we get
\begin{align}\label{a-density-constants}
\tilde c_k
=& c_k^A exp\Bigl(-k \frac{N(N-1)}{2}(1+\ln 2)+k\sum_{j=1}^N j\ln j\Bigr)\cdot (2k)^{kN(N-1)/2}
\notag\\
=& c_k^A  (k/e)^{kN(N-1)/2}\cdot\prod_{j=1}^N j^{kj}\notag\\
=&\frac{N!}{(2\pi)^{N/2}} \cdot\prod_{j=1}^{N}\frac{\Gamma(1+k)}{\Gamma(1+jk)}\cdot
   (k/e)^{kN(N-1)/2}\cdot\prod_{j=1}^N j^{kj}\notag\\
=&\frac{1}{(2\pi)^{N/2}} \cdot\prod_{j=1}^{N}\frac{\Gamma(k)}{\Gamma(jk)}\cdot
   (k/e)^{kN(N-1)/2}\cdot\prod_{j=1}^N j^{kj}.
\end{align}
 Stirling's series (see e.g.~Section 12.33 of \cite{WW}) now implies  that for $x>1$,
\begin{equation}\label{stirling}
  \Gamma(x) = \sqrt{2\pi}\,x^{x-1/2}\, e^{-x}\, e^{\mu(x)}  \quad\quad\text{with}  \quad\quad
  \frac{1}{ 12\>  x} - \frac{1}{ 360\>  x^3} < \mu( x ) < \frac{1}{ 12\>  x}.
  \end{equation}
We conclude with elementary calculus that
\begin{equation}\label{a-density-constants-limit}
\tilde c_k=\frac{\sqrt{N!}}{(2\pi)^{N/2}}\cdot 
e^{M}  \quad\quad\text{with}  \quad\quad  M\ge \frac{N-1}{26k}.
\end{equation}
In fact, for the exponent $M$ we obtain from (\ref{a-density-constants}) and (\ref{stirling}) that for $k\ge1$,
\begin{align}
  M&= (N-1)\mu(k)-\sum_{l=2}^N \mu(lk)\notag\\
  &\ge \frac{N-1}{12k}- \frac{1}{12}\sum_{l=2}^N \frac{1}{lk}- \frac{N-1}{360k^3}=\frac{1}{12k}\sum_{l=2}^N (1-\frac{1}{l})- \frac{N-1}{360k^3}
\notag\\
&\ge\frac{N-1}{24k}- \frac{N-1}{360k}=\frac{14(N-1)}{360k}\ge\frac{N-1}{26k}.\notag
\end{align}

We next turn to the function  $h_k(y)$. The Lagrange remainder  of the power series of $\ln(1+x)$ for
  $|x|\le 1/2$ yields that
$$\ln(1+x)= x-x^2/2+x^3/3- 4r(x)x^4\quad\quad\text{with}  \quad\quad r(x)\in[0,1].$$
We now apply this to 
\begin{equation}\label{special-x-a}x=\frac{y_i-y_j}{\sqrt{2k}(z_{i,N}-z_{j,N})} \end{equation}
for $i<j$ with $|x|\le1/2$,   and we use that, by Lemma \ref{char-zero-A}, 
$$%\begin{equation}\label{zero-equation-a}
-\sqrt{2k} \langle y,{\bf z}_N\rangle+\sqrt{2k} \sum_{i<j} \frac{y_i-y_j}{z_{i,N}-z_{j,N}}=
\sqrt{2k}\sum_{i=1}^N y_i\bigl(-z_i+ \sum_{j: \> j\ne i} \frac{1}{z_{i,N}-z_{j,N}}\bigr)=0.$$
%\end{equation}
We conclude that
\begin{equation}\label{density-a-limit-1}
h_k(y)= exp\Bigl(-\|y\|_2^2/2-\frac{1}{2}\sum_{i<j}\frac{(y_i-y_j)^2}{(z_{i,N}-z_{j,N})^2} + \frac{1}{k^{1/2}}R_1(y) -\frac{1}{k}R_2(y)\Bigr)
\end{equation}
with 
\begin{equation}\label{r1-r2-a}
R_1(y):=\frac{1}{3\sqrt 2}\sum_{i<j}\frac{(y_i-y_j)^3}{(z_{i,N}-z_{j,N})^3} \quad\text{and}\quad
0\le R_2(y)\le2\sum_{i<j}\frac{(y_i-y_j)^4}{(z_{i,N}-z_{j,N})^4}.
\end{equation}
In particular, we see from (\ref{sum-inverse-herm-2})  that for $y\in\mathbb R^N$ with sup-norm $\|y\|_\infty\le\epsilon$,
\begin{equation}\label{density-a-limit-2}
  0\le R_2(y)\le   (2   \|y\|_\infty)^4 \frac{N^3}{6}\le  \frac{8}{3} \epsilon^4 N^3.
  \end{equation}
Notice here that in (\ref{r1-r2-a}) only one half of the summands of the left hand side in  (\ref{sum-inverse-herm-2}) appears.
We here also notice that for $\|y\|_\infty\le\epsilon$ and $\epsilon$ with
\begin{equation}\label{condition-epsilon-a}0<\epsilon\le \frac{\sqrt{2k}}{2\sqrt N},\end{equation}
 Lemma \ref{min-distance-a}
ensures
that the computations above for $x$ as in 
(\ref{special-x-a}) can be applied.

We now compare the density $f_k^A(y)= \tilde c_k\cdot h_k(y)$ with 
the normal distribution $N(0,\Sigma_N)$
in Theorem \ref{clt-main-a} with the density
$$\frac{\sqrt{N!}}{(2\pi)^{N/2}} exp\Bigl(-\|y\|_2^2/2-\frac{1}{2}\sum_{i<j}\frac{(y_i-y_j)^2}{(z_{i,N}-z_{j,N})^2}\Bigr).$$
In particular, (\ref{density-a-limit-1}), (\ref{density-a-limit-2}) and (\ref{a-density-constants-limit}) show that
\begin{align}\label{main-estimate-a1}
  P\Bigl(\|X_{k,N}-\sqrt{2k}{\bf z}_N\|_2\le\epsilon\Bigr)&=\int_{\{y\in\mathbb R^N:\> \|y\|_2\le\epsilon\}} f_k^A(y)\> dy
  \\
  &\ge e^{\frac{N-1}{30k}} \int_{\{y\in\mathbb R^N:\> \|y\|_2\le\epsilon\}}e^{R_1(y)/\sqrt k}e^{-R_2(y)/k} \> dN(0,\Sigma_N)(y)
  \notag\\
  &\ge \tilde E(N,\epsilon,k)\int_{\{y\in\mathbb R^N:\> \|y\|_2\le\epsilon\}}e^{R_1(y)/\sqrt k} \> dN(0,\Sigma_N)(y)
  \notag
\end{align}
with
\begin{equation}\label{E-tilde}
\tilde E(N,\epsilon,k):= 
exp\Bigl(\frac{1}{k}\Bigl(\frac{N-1}{26}- \frac{8}{3} \epsilon^4N^3\Bigr)\Bigr).
\end{equation}
As the density of  $N(0,\Sigma_N)$ and the integration domain are
invariant under the transform $y\mapsto -y$ on the right hand side of (\ref{main-estimate-a1}),
the integral on the right hand side  of (\ref{main-estimate-a1}) is also equal to
$$  \int_{\{y\in\mathbb R^N:\> \|y\|_2\le\epsilon\}}e^{-R_1(y)/k^{1/2}} \> dN(0,\Sigma_N)(y)$$
and thus to
$$  \int_{\{y\in\mathbb R^N:\> \|y\|_2\le\epsilon\}}\cosh(R_1(y)/k^{1/2})\> dN(0,\Sigma_N)(y).$$
 We thus conclude from (\ref{main-estimate-a1}) and $\cosh u\ge1$ for $u\in\mathbb R$ that 
\begin{equation}\label{main-estimate-a2}
  P\Bigl(\|X_{k,N}-\sqrt{2k}{\bf z}_N\|_2\le\epsilon\Bigr)\ge
 \tilde E(N,\epsilon,k)\cdot N(0,\Sigma_N)(\{y\in\mathbb R^N:\> \|y\|_2\le\epsilon\}).
\end{equation}
Therefore, 
%for  $\epsilon>0$ sufficiently small with $ E(N,\epsilon,k)\ge1$,
\begin{align}\label{main-estimate-a3}
  P\Bigl(\|X_{k,N}-\sqrt{2k}{\bf z}_N\|_2>\epsilon\Bigr)&=1- P\Bigl(\|X_{k,N}-z\|_2\le\epsilon\Bigr)\notag\\
  &\le 1- \tilde E(N,\epsilon,k)\cdot N(0,\Sigma_N)(\{y\in\mathbb R^N:\> \|y\|_2\le\epsilon\})
  \notag\\
&= 1- \tilde E(N,\epsilon,k)\Bigl(1- N(0,\Sigma_N)(\{y\in\mathbb R^N:\> \|y\|_2>\epsilon\})\Bigr)\notag\\
%&\le E(N,\epsilon,k)\cdot N(0,\Sigma_N)(\{y\in\mathbb R^N:\> \|y\|_2>\epsilon\})\Bigr)\notag\\
&= 1- \tilde E(N,\epsilon,k)\Bigl(1- P(\|\tilde X\|_2>\epsilon)\Bigr)
\end{align}
with some $N(0, diag(1,1/2,\ldots,1/N))$-distributed random variable $\tilde X$
due to the eigenvalues of $\Sigma_N$ in Theorem \ref{clt-main-a}. We now switch to a $N(0,I_N)$-distributed
random variable $\hat X=(\hat X_1,\ldots,\hat X_N)$ and obtain from
the Markov inequality that
\begin{align}\label{main-estimate-a4}
  P(\|\tilde X\|_2>\epsilon)&=
P\Bigl(\sum_{l=1}^N \hat X_l^2/l>\epsilon^2\Bigr)
\le e^{-\delta\epsilon^2} \cdot\mathbb E\Bigl(exp\Bigl(\delta\sum_{l=1}^N \hat X_l^2/l\Bigr)\Bigr)\notag\\
&=e^{-\delta\epsilon^2}\cdot \prod_{l=1}^N \mathbb E\Bigl(exp(\delta \hat X_l^2/l)\Bigr)
\end{align}
for $\delta>0$. The right hand side of this estimation is finite  for $2\delta<1$ with
$\mathbb E(exp(\delta \hat X_l^2/l))=(1-2\delta/l)^{-1/2}$ by an elementary calculation.
We next estimate the right hand side of  (\ref{main-estimate-a4}). For this we use
$$ \sum_{l=2}^N \frac{1}{l}\le \ln N,\quad\quad  \sum_{l=1}^N \frac{1}{l^2}\le  2$$
and obtain 
from the power series of $\ln(1-x)$  for $2\delta<1$ that
\begin{align}
  \prod_{l=1}^N(1-2\delta/l)^{-1/2}&= exp\Bigl(-\frac{1}{2}\sum_{l=1}^N\ln(1-2\delta/l)\Bigr)
  =exp\Bigl(\frac{1}{2}\sum_{l=1}^N\sum_{m=1}^\infty\frac{(2\delta)^m}{l^m\cdot m}\Bigr)\notag\\
  &\le exp\Bigl((1+\ln N)\delta +\frac{1}{2}\sum_{m=2}^\infty\Bigl(\frac{(2\delta)^m}{m}\sum_{l=1}^\infty l^{-m}
  \Bigr) \Bigr)\notag\\
  &\le   exp\Bigl((1+\ln N)\delta +\sum_{m=2}^\infty\frac{(2\delta)^m}{m} \Bigr)\notag\\
  &=  exp\Bigl((1+\ln N)\delta -2\delta-\ln(1-2\delta)\Bigr)\notag
\end{align}
Hence, by (\ref{main-estimate-a4}),
\begin{equation}\label{main-estimate-a5}
  P(\|\tilde X\|_2>\epsilon)\le
 \frac{1}{1-2\delta}  exp\Bigl((-\epsilon^2-1+\ln N)\delta\Bigr).
\end{equation}
In combination with (\ref{main-estimate-a3}) we obtain that
\begin{align}\label{main-estimate-a6}
 P\Bigl(\|X_{k,N}-\sqrt{2k}\cdot {\bf z}_N\|_2>\epsilon\Bigr)
 \le& 1- \tilde E(N,\epsilon,k) \\&+ 
 \frac{\tilde  E(N,\epsilon,k)}{1-2\delta}  \cdot exp\Bigl((-\epsilon^2-1+\ln N)\delta\Bigr).
 \notag \end{align}

We now replace $\epsilon$ in (\ref{main-estimate-a6}) by $\sqrt{2k}\cdot\epsilon$. Then
the condition (\ref{main-condition-prp-a})
in the proposition yields the condition  (\ref{condition-epsilon-a}), and
$E(N,\epsilon,k):=\tilde E(N,\sqrt{2k}\cdot\epsilon,k)$ and (\ref{main-estimate-a6}) lead to 
\begin{equation}\label{main-estimate-a7}
 P\Bigl(\|\frac{X_{k,N}}{\sqrt{2k}}-{\bf z}_N\|_2>\epsilon\Bigr)
 \le 1- E(N,\epsilon,k) 
 + 
 \frac{E(N,\epsilon,k)}{1-2\delta} \cdot e^{(-2k\epsilon^2-1+\ln N)\delta}.
   \end{equation}
We now optimize (\ref{main-estimate-a7}) w.r.t. $\delta\in[0,1[/2[$ and study the function
   \begin{equation}\label{main-estimate-a8}\phi(\delta):= \frac{1}{1-2\delta} \cdot e^{r\delta} \quad\quad\text{with}\quad\quad
        r:=-2k\epsilon^2-1+\ln N .\end{equation}
   Elementary calculus shows that for $r\ge -2$, the function
   $\phi$ is minimal for $\delta=0$, in which case the r.h.s.~of (\ref{main-estimate-a7})
   is equal to 1,
   and thus (\ref{main-estimate-a7}) is not interesting in this case.
   On the other hand, for $r<-2$, $\phi$ has its minimum on $[0,1[/2[$ at $\delta=\frac{1}{2}+\frac{1}{r}$.
         We now insert this into (\ref{main-estimate-a7}) and apply $1-e^{-x}\le x$ to
         $x:=\frac{32}{3}\epsilon^4kN^3-\frac{N-1}{26\cdot k}\in\mathbb R$ such that (\ref{main-estimate-a7})
         leads to  the claim (\ref{main-estimate-prp-a}).
\end{proof}

\end{document}